\documentclass{amsart}

\usepackage{amsmath,amsfonts,amssymb,amsmath,latexsym,multirow,color}

\newtheorem{theorem}{Theorem}[section]

\theoremstyle{definition}

\theoremstyle{remark}

\numberwithin{equation}{section}

\begin{document}

\title[Symmetry identities for generalized twisted Euler polynomials
]{$\begin{array}{c}
         \text{Symmetry identities for generalized twisted Euler polynomials}\\
           \text{twisted by ramified roots of unity}
       \end{array}$
 }

\author{dae san kim}
\address{Department of Mathematics, Sogang University, Seoul 121-742, Korea}
\curraddr{Department of Mathematics, Sogang University, Seoul
121-742, Korea} \email{dskim@sogong.ac.kr}

\begin{abstract}
We derive eight identities of symmetry in three variables related to
generalized twisted Euler polynomials and alternating generalized
twisted power sums, both of which are twisted by ramified roots of
unity.  All of these are new, since there have been results only
about identities of symmetry in two variables. The derivations of
identities are based on the $p$-adic integral expression of the
generating function for the generalized twisted Euler polynomials
and the quotient of $p$-adic integrals that can be expressed as the
exponential generating function for the alternating generalized
twisted power sums.

\end{abstract}


\date{}

\dedicatory{ }


\maketitle

\section{Introduction and preliminaries}
Let $p$ be a fixed odd prime. Throughout this paper, $\mathbb{Z}_p$,
$\mathbb{Q}_p$, $\mathbb{C}_p$ will respectively denote the ring of
$p$-adic integers, the field of $p$-adic rational numbers and the
completion of the algebraic closure of $\mathbb{Q}_p$. Assume that
$|~|_p$ is the normalized absolute value of $\mathbb{C}_p$, such
that $|p|_p=\frac{1}{p}$. The group $\Gamma$ of all roots of unity
of $\mathbb{C}_p$ is the direct product of its subgroups $\Gamma_u$
(the subgroup of unramified roots of unity) and $\Gamma_r$ (the
subgroup of ramified roots of unity). Namely,

\begin{equation*}
\Gamma=\Gamma_u\cdot\Gamma_r,~\Gamma_u\cap\Gamma_r=\{1\},
\end{equation*}
where
\begin{equation*}
\Gamma_u=\{\xi\in\mathbb{C}_p~|~\xi^s=1~for~
some~s\in\mathbb{Z}_{>0}~with~(s,p)=1\},
\end{equation*}

\begin{equation*}
\begin{split}
&\Gamma_r=\{\xi\in\mathbb{C}_p~|~\xi^{p^s}=1~for~
some~s\in\mathbb{Z}_{>0}\}.\\
&
\end{split}
\end{equation*}

\noindent Let $d$ be a fixed odd positive integer. Then we let
\begin{equation*}
\begin{split}
&\\
&X=X_d=\lim_{\overleftarrow{N}}\mathbb{Z}/dp^N \mathbb{Z}
=\cup^{dp-1}_{a=0}(a+dp\mathbb{Z}_p),
\end{split}
\end{equation*}
with
\begin{equation*}
\begin{split}
&a+dp^N\mathbb{Z}_p=\{x\in X~|~x\equiv a~(mod~dp^N)\},\\
&
\end{split}
\end{equation*}
and let $\pi:X\rightarrow\mathbb{Z}_p$ be the map given by the
inverse limit of the natural maps
\begin{equation*}
\begin{split}
&\\
&\mathbb{Z}/dp^N\mathbb{Z}~\rightarrow~\mathbb{Z}/p^N\mathbb{Z}.\\
&
\end{split}
\end{equation*}

\noindent If $g$ is a function on $\mathbb{Z}_p$, we will use the
same notation to denote the function $g\circ\pi$. Let
$\chi:(\mathbb{Z}/d\mathbb{Z})^*\rightarrow\overline{Q}^*$ be a
(primitive) Dirichlet character of conductor $d$. Then it will be
pulled back to $X$ via the natural map
$X\rightarrow\mathbb{Z}/d\mathbb{Z}$. Here we fix, once and for all,
an imbedding $\overline{Q}\rightarrow\mathbb{C}_p$, so that $\chi$
is regarded as a map of $X$ to $\mathbb{C}_p$.(cf. [8]).

For a continuous function $f:X\rightarrow\mathbb{C}_p$, the $p$-adic
fermionic integral of $f$ is defined by
\begin{equation*}
\int_{X}f(z)d\mu_{-1}(z)=\lim_{N\rightarrow\infty}\sum_{j=0}^{dp^N-1}f(j)(-1)^j.
\end{equation*}

\noindent Then it is easy to see that

\begin{equation}\label{N1}
\int_{X}f(z+1)d\mu_{-1}(z)+\int_{X}f(z)d\mu_{-1}(z)=2f(0).
\end{equation}

\noindent More generally, we deduce from (\ref{N1}) that, for any
odd positive integer $n$,

\begin{equation}\label{N2}
\int_{X}f(z+n)d\mu_{-1}(z)+\int_{X}f(z)d\mu_{-1}(z)=2\sum_{a=0}^{n-1}
(-1)^af(a),
\end{equation}

\noindent and that, for any even positive integer $n$,

\begin{equation*}
\int_{X}f(z+n)d\mu_{-1}(z)-\int_{X}f(z)d\mu_{-1}(z)=2\sum_{a=0}^{n-1}
(-1)^{a-1}f(a).
\end{equation*}

\noindent  Throughout this paper, we let $\xi\in\Gamma_r$ be any
fixed root of unity, and let

\begin{equation}\label{N3}
E=\left\{ t\in \mathbb{C}_p ~\Big|~
|t|_p<p^{-\frac{1}{p-1}}\right\}.
\end{equation}

\noindent Then, for each fixed $t\in E$, the function $e^{zt}$ is
analytic on $\mathbb{Z}_p$ and hence considered as a function on X,
and, by applying (\ref{N2}) to $f$ with $f(z)=\chi(z)\xi^ze^{zt}$,
we get the $p$-adic integral expression of the generating function
for the generalized twisted Euler numbers $E_{n,\chi,\xi}$ attached
to $\chi$ and $\xi$:

\begin{equation}\label{N4}
\int_{X}\chi(z)\xi^ze^{zt}d\mu_{-1}(z)=\frac{2}{\xi^de^{dt}+1}\sum_{a=0}^{d-1}(-1)^a
\chi(a)\xi^ae^{at}=\sum_{n=0}^{\infty}E_{n,\chi,\xi}\frac{t^n}{n!}
~(t\in E).
\end{equation}

\noindent So we have the following $p$-adic integral expression of
the generating function for the generalized twisted Euler
polynomials $E_{n,\chi}(x)$ attached to $\chi$ and $\xi$:

\begin{equation}\label{N5}
\begin{split}
\int_{X}\chi(z)&\xi^ze^{(x+z)t}d\mu_{-1}(z)\\
&=\frac{2e^{xt}}{\xi^de^{dt}+1}\sum_{a=0}^{d-1}(-1)^a
\chi(a)\xi^ae^{at}=\sum_{n=0}^{\infty}E_{n,\chi,\xi}(x)\frac{t^n}{n!}
~(t\in E,~x\in \mathbb{Z}_p).
\end{split}
\end{equation}

\noindent Also, from (\ref{N1}) we have:

\begin{equation}\label{N6}
\int_{X}\xi^ze^{zt}d\mu_{-1}(z)=\frac{2}{\xi e^{t}+1}~(t\in E).
\end{equation}

Let $T_k(n;\chi,\xi)$ denote the $k$th alternating generalized
twisted power sum of the first $n+1$ nonnegative integers attached
to $\chi$ and $\xi$, namely

\begin{equation}\label{N7}
T_k(n;\chi,\xi)=\sum_{a=0}^{n}(-1)^a\chi(a)\xi^aa^k.
\end{equation}

From (\ref{N4}), (\ref{N6}), and (\ref{N7}), one easily derives the
following identities: for any odd positive integer $w$,

\begin{equation}\label{N8}
\frac
  {\int_{X}\chi(x)\xi^xe^{xt}d\mu_{-1}(x)}{\int_{X}\xi^{wdy}e^{wdyt}d\mu_{-1}(y)}
=\frac{\xi^{wd}e^{wdt}+1}{\xi^de^{dt}+1}\sum_{a=0}^{d-1}(-1)^a\chi(a)\xi^ae^{at}
\end{equation}

\begin{equation}\label{N9}
\quad\qquad\qquad=\sum_{a=0}^{wd-1}(-1)^a \chi(a)\xi^ae^{at}
\end{equation}

\begin{equation}\label{N10}
\qquad\qquad\qquad\qquad\qquad=\sum_{k=0}^{\infty}T_k(wd-1;\chi,\xi)\frac{t^k}{k!}~(t\in
E).
\end{equation}

In what follows, we will always assume that the $p$-adic integrals
of the various (twisted) exponential functions on $X$ are defined
for $t\in E$(cf. (\ref{N3})), and therefore it will not be
mentioned.

\cite{B1}, \cite{B2}, \cite{B6}, \cite{B9} and \cite{B10} are some
of the previous works on identities of symmetry involving Bernoulli
polynomials and power sums. For the brief history, one is referred
to those papers. For the first time, the idea of \cite{B6} was
extended in [4, 5] to  the case three variables so as to yield many
new identities with abundant symmetry. This added some new
identities of symmetry even to the existing ones in two variables as
well.

In this paper, we will produce  8 identities of symmetry in three
variables $w_1$, $w_2$, $w_3$ related to generalized twisted Euler
polynomials and alternating generalized twisted power sums, both of
which are twisted by ramified roots of unity(i.e., $p$-power roots
of unity)(cf. (\ref{N39})-(\ref{N42}), (\ref{N45})-(\ref{N48})). All
of these seem to be new, since there have been results only about
identities of symmetry in two variables in the
literature(\cite{B7}). On the other hand, in \cite{B3} the measure
introduced by Koblitz(cf. \cite{B8})was modified in order to treat
the unramified roots of unity case (i.e., the orders of the roots of
unity are prime to $p$ and the conductors of Dirichlet characters).

The following is stated as Theorem \ref{T2} and  an example of the
full six symmetries in $w_1$, $w_2$, $w_3$.

\begin{align*}
&\sum_{k+l+m=n}{\binom{n}{k,l,m}}E_{k,\chi,\xi^{w_2w_3}}(w_1y_1)
E_{l,\chi,\xi^{w_1w_3}}(w_2y_2)\\
&\qquad\qquad\qquad\qquad\qquad\qquad\qquad
\times T_{m}(w_3d-1;\chi,\xi^{w_1w_2})w_{1}^{l+m}w_{2}^{k+m}w_{3}^{k+l}\\
&\sum_{k+l+m=n}{\binom{n}{k,l,m}}E_{k,\chi,\xi^{w_2w_3}}(w_1y_1)
E_{l,\chi,\xi^{w_1w_2}}(w_3y_2)\\
&\qquad\qquad\qquad\qquad\qquad\qquad\qquad
\times T_{m}(w_2d-1;\chi,\xi^{w_1w_3})w_{1}^{l+m}w_{3}^{k+m}w_{2}^{k+l}\\
&\sum_{k+l+m=n}{\binom{n}{k,l,m}}E_{k,\chi,\xi^{w_1w_3}}(w_2y_1)
E_{l,\chi,\xi^{w_2w_3}}(w_1y_2)\\
&\qquad\qquad\qquad\qquad\qquad\qquad\qquad
\times T_{m}(w_3d-1;\chi,\xi^{w_1w_2})w_{2}^{l+m}w_{1}^{k+m}w_{3}^{k+l}\\
&\sum_{k+l+m=n}{\binom{n}{k,l,m}}E_{k,\chi,\xi^{w_1w_3}}(w_2y_1)
E_{l,\chi,\xi^{w_1w_2}}(w_3y_2)\\
&\qquad\qquad\qquad\qquad\qquad\qquad\qquad
\times T_{m}(w_1d-1;\chi,\xi^{w_2w_3})w_{2}^{l+m}w_{3}^{k+m}w_{1}^{k+l}\\
&\sum_{k+l+m=n}{\binom{n}{k,l,m}}E_{k,\chi,\xi^{w_1w_2}}(w_3y_1)
E_{l,\chi,\xi^{w_1w_3}}(w_2y_2)\\
&\qquad\qquad\qquad\qquad\qquad\qquad\qquad
\times T_{m}(w_1d-1;\chi,\xi^{w_2w_3})w_{3}^{l+m}w_{2}^{k+m}w_{1}^{k+l}\\
&\sum_{k+l+m=n}{\binom{n}{k,l,m}}E_{k,\chi,\xi^{w_1w_2}}(w_3y_1)
E_{l,\chi,\xi^{w_2w_3}}(w_1y_2)\\
&\qquad\qquad\qquad\qquad\qquad\qquad\qquad
\times T_{m}(w_2d-1;\chi,\xi^{w_1w_3})w_{3}^{l+m}w_{1}^{k+m}w_{2}^{k+l}.\\
\end{align*}

The derivations of identities are based on the $p$-adic integral
expression of the generating function for the generalized twisted
Euler polynomials in (\ref{N5}) and the quotient of integrals in
(\ref{N8})-(\ref{N10}) that can be expressed as the exponential
generating function for the alternating generalized twisted power
sums. These abundance of symmetries would not be unearthed if such
$p$-adic integral representations had not been available. We
indebted this idea to the papers [6, 7].

\section{Several types of quotients of $p$-adic fermionic integrals}
Here we will introduce several types of quotients of $p$-adic
fermionic integrals on $X$ or $X^3$  from which some interesting
identities follow owing to the built-in symmetries in $w_1$, $w_2$,
$w_3$. In the following, $w_1$, $w_2$, $w_3$ are all positive
integers and all of the explicit expressions of integrals in
(\ref{N12}), (\ref{N14}), (\ref{N16}), and (\ref{N18}) are obtained
from the identities in (\ref{N4}) and (\ref{N6}). To ease notations,
from now on we will suppress $\mu_{-1}$ and denote, for example,
$d\mu_{-1}(x)$ and $d\mu_{-1}(x_1)d\mu_{-1}(x_2)d\mu_{-1}(x_3)$
respectively simply by $dx$ and $dX$.

\noindent (a) Type $\Lambda_{23}^{i}$ (for $i=0,1,2,3$)
\begin{equation}\label{N11}
\begin{split}
&I(\Lambda_{23}^{i})\\
&=\frac{\int_{X^3}\chi(x_1)\chi(x_2)\chi(x_3)\xi^{w_2w_3x_1+w_1w_3x_2+w_1w_2x_3}
e^{(w_2w_3x_1+w_1w_3x_2+w_1w_2x_3+w_1w_2w_3(\sum_{j=1}^{3-i}y_j))t}dX}
 {(\int_{X}\xi^{dw_1w_2w_3x_4}e^{dw_1w_2w_3x_4t}dx_4)^i}\\
 \end{split}
\end{equation}
\begin{equation}\label{N12}
\begin{split}
&=\frac{2^{3-i}e^{w_1w_2w_3(\sum_{j=1}^{3-i}y_j)t}(\xi^{dw_1w_2w_3}e^{dw_1w_2w_3t}+1)^i}
{(\xi^{dw_2w_3}e^{dw_2w_3t}+1)(\xi^{dw_1w_3}e^{dw_1w_3t}+1)(\xi^{dw_1w_2}e^{dw_1w_2t}+1)}\\
&\qquad\qquad\times(\sum_{a=0}^{d-1}(-1)^a\chi(a)\xi^{aw_2w_3}e^{aw_2w_3t})(\sum_{a=0}^{d-1}(-1)^a\chi(a)\xi^{aw_1w_3}e^{aw_1w_3t})\\
&\qquad\qquad\qquad\qquad\qquad\qquad\qquad\qquad\qquad\times(\sum_{a=0}^{d-1}(-1)^a\chi(a)\xi^{aw_1w_2}e^{aw_1w_2t}).
\end{split}
\end{equation}

\noindent (b) Type $\Lambda_{13}^{i}$ (for $i=0,1,2,3$)
\begin{equation}\label{N13}
I(\Lambda_{13}^{i})=\frac{\int_{X^3}\chi(x_1)\chi(x_2)\chi(x_3)
\xi^{w_1x_1+w_2x_2+w_3x_3}
e^{(w_1x_1+w_2x_2+w_3x_3+w_1w_2w_3(\sum_{j=1}^{3-i}y_j))t}dX}
 {(\int_{X}\xi^{dw_1w_2w_3x_4}e^{dw_1w_2w_3x_4t}dx_4)^i}\\
\end{equation}
\begin{equation}\label{N14}
\begin{split}
&=\frac{2^{3-i}e^{w_1w_2w_3(\sum_{j=1}^{3-i}y_j)t}(\xi^{dw_1w_2w_3}e^{dw_1w_2w_3t}+1)^i}
{(\xi^{dw_1}e^{dw_1t}+1)(\xi^{dw_2}e^{dw_2t}+1)(\xi^{dw_3}e^{dw_3t}+1)}\\
&\times(\sum_{a=0}^{d-1}(-1)^a\chi(a)\xi^{aw_1}e^{aw_1t})(\sum_{a=0}^{d-1}(-1)^a\chi(a)\xi^{aw_2}e^{aw_2t})
(\sum_{a=0}^{d-1}(-1)^a\chi(a)\xi^{aw_3}e^{aw_3t}).
\end{split}
\end{equation}

\noindent (c-0) Type $\Lambda_{12}^{0}$
\begin{equation}\label{N15}
\begin{split}
I&(\Lambda_{12}^{0})\\
&=\int_{X^3}\chi(x_1)\chi(x_2)\chi(x_3) \xi^{w_1x_1+w_2w_2+w_3x_3}
e^{(w_1x_1+w_2x_2+w_3w_3+w_2w_3y+w_1w_3y+w_1w_2y)t}dX
\end{split}
\end{equation}
\begin{equation}\label{N16}
\begin{split}
&=\frac{8e^{(w_2w_3+w_1w_3+w_1w_2)yt}}
{(\xi^{dw_1}e^{dw_1t}+1)(\xi^{dw_2}e^{dw_2t}+1)(\xi^{dw_3}e^{dw_3t}+1)}\\
&\times(\sum_{a=0}^{d-1}(-1)^a\chi(a)\xi^{aw_1}e^{aw_1t})(\sum_{a=0}^{d-1}(-1)^a\chi(a)\xi^{aw_2}e^{aw_2t})
(\sum_{a=0}^{d-1}(-1)^a\chi(a)\xi^{aw_3}e^{aw_3t}).
\end{split}
\end{equation}

\noindent (c-1) Type $\Lambda_{12}^{1}$
\begin{equation}\label{N17}
I(\Lambda_{12}^{1})=\frac{\int_{X^3}\chi(x_1)\chi(x_2)\chi(x_3)
\xi^{w_1x_1+w_2x_2+w_3x_3} e^{(w_1x_1+w_2x_2+w_3x_3)t}dx_1dx_2dx_3}
 {\int_{X^3}\xi^{d(w_2w_3z_1+w_1w_3z_2+w_1w_2z_3)}
 e^{d(w_2w_3z_1+w_1w_3z_2+w_1w_2z_3)t}dz_1dz_2dz_3}
\end{equation}
\begin{equation}\label{N18}
\begin{split}
&=\frac{(\xi^{dw_2w_3}e^{dw_2w_3t}+1)(\xi^{dw_1w_3}e^{dw_1w_3t}+1)(\xi^{dw_1w_2}e^{dw_1w_2t}+1)}
{(\xi^{dw_1}e^{dw_1t}+1)(\xi^{dw_2}e^{dw_2t}+1)(\xi^{dw_3}e^{dw_3t}+1)}\\
&\times(\sum_{a=0}^{d-1}(-1)^a\chi(a)\xi^{aw_1}e^{aw_1t})
(\sum_{a=0}^{d-1}(-1)^a\chi(a)\xi^{aw_2}e^{aw_2t})
(\sum_{a=0}^{d-1}(-1)^a\chi(a)\xi^{aw_3}e^{aw_3t}).
\end{split}
\end{equation}

All of the above $p$-adic integrals of various types are invariant
under all permutations of $w_1$, $w_2$, $w_3$, as one can see either
from $p$-adic integral representations in (\ref{N11}), (\ref{N13}),
(\ref{N15}), and (\ref{N17}) or from their explicit evaluations in
(\ref{N12}), (\ref{N14}), (\ref{N16}), and (\ref{N18}).

\section{Identities for generalized twisted Euler polynomials}

In the following $w_1$, $w_2$, $w_3$ are all odd positive integers
except for (a-0) and (c-0), where they are any positive integers.
First, let's consider Type $\Lambda_{23}^{i}$, for each $i=0,1,2,3$.
The following results can be easily obtained from (\ref{N5}) and
(\ref{N8})-(\ref{N10}).

\noindent (a-0)
\begin{equation}\label{N19}
\begin{split}
I(&\Lambda_{23}^{0})\\
&=\int_{X}\chi(x_1)\xi^{w_2w_3x_1}e^{w_2w_3(x_1+w_1y_1)t}dx_1
\int_{X}\chi(x_2)\xi^{w_1w_3x_2}e^{w_1w_3(x_2+w_2y_2)t}dx_2\\
&\qquad\qquad\qquad\qquad\qquad\qquad\qquad\qquad\times
\int_{X}\chi(x_3)\xi^{w_1w_2x_3}e^{w_1w_2(x_3+w_3y_3)t}dx_3\\
&=(\sum_{k=0}^{\infty}\frac{E_{k,\chi,\xi^{w_2w_3}}(w_1y_1)}{k!}(w_2w_3t)^k)
(\sum_{l=0}^{\infty}\frac{E_{l,\chi,\xi^{w_1w_3}}(w_2y_2)}{l!}(w_1w_3t)^l)\\
&\qquad\qquad\qquad\qquad\qquad\qquad\qquad\qquad\times
(\sum_{m=0}^{\infty}\frac{E_{m,\chi,\xi^{w_1w_2}}(w_3y_3)}{m!}(w_1w_2t)^m)\\
&=\sum_{n=0}^{\infty}(\sum_{k+l+m=n}\binom{n}{k,l,m}E_{k,\chi,\xi^{w_2w_3}}(w_1y_1)
E_{l,\chi,\xi^{w_1w_3}}(w_2y_2)\\
&\qquad\qquad\qquad\qquad\qquad\qquad\qquad\times
E_{m,\chi,\xi^{w_1w_2}}(w_3y_3)
w_{1}^{l+m}w_{2}^{k+m}w_{3}^{k+l})\frac{t^n}{n!},
\end{split}
\end{equation}

\noindent where the inner sum is over all nonnegative integers
$k,~l,~m$ with $k+l+m=n$, and
\begin{equation}\label{N20}
\binom{n}{k,l,m}=\frac{n!}{k!l!m!}.
\end{equation}

\noindent (a-1) Here we write $I(\Lambda_{23}^{1})$ in two different
ways:

\begin{equation}\label{N21}
\begin{split}
\noindent
(1)~I(\Lambda_{23}^{1})&=\int_{X}\chi(x_1)\xi^{w_2w_3x_1}e^{w_2w_3(x_1+w_1y_1)t}dx_1\\
&\qquad\qquad\qquad\qquad \times
\int_{X}\chi(x_2)\xi^{w_1w_3x_2}e^{w_1w_3(x_2+w_2y_2)t}dx_2\\
&\qquad\qquad\qquad\qquad\qquad\qquad\qquad
\times\frac{\int_{X}\chi(x_3)\xi^{w_1w_2x_3}e^{w_1w_2x_3t}dx_3}{\int_{X}\xi^{dw_1w_2w_3x_4}e^{dw_1w_2w_3x_4t}dx_4}\\
&=(\sum_{k=0}^{\infty}E_{k,\chi,\xi^{w_2w_3}}(w_1y_1)\frac{(w_2w_3t)^k}{k!})
(\sum_{l=0}^{\infty}E_{l,\chi,\xi^{w_1w_3}}(w_2y_2)\frac{(w_1w_3t)^l}{l!})\\
&\qquad\qquad\qquad\qquad\qquad
\times(\sum_{m=0}^{\infty}T_m(w_3d-1;\chi,\xi^{w_1w_2})\frac{(w_1w_2t)^m}{m!})\\
\end{split}
\end{equation}

\begin{equation}\label{N22}
\begin{split}
\qquad\qquad\qquad=&\sum_{n=0}^{\infty}(\sum_{k+l+m=n}\binom{n}{k,l,m}E_{k,\chi,\xi^{w_2w_3}}(w_1y_1)E_{l,\chi,\xi^{w_1w_3}}(w_2y_2)\\
&\qquad\qquad\qquad\quad\times
T_m(w_3d-1;\chi,\xi^{w_1w_2})w_{1}^{l+m}w_{2}^{k+m}w_{3}^{k+l})\frac{t^n}{n!}.
\end{split}
\end{equation}

\noindent (2) Invoking (\ref{N9}), (\ref{N21}) can also be written
as

\begin{equation}\label{N23}
\begin{split}
I(&\Lambda_{23}^{1})\\
&=\sum_{a=0}^{w_3d-1}(-1)^a\chi(a)\xi^{aw_1w_2}
\int_{X}\chi(x_1)\xi^{w_2w_3x_1}e^{w_2w_3(x_1+w_1y_1)t}dx_1\\
&\qquad\qquad\qquad\qquad\qquad\times
\int_{X}\chi(x_2)\xi^{w_1w_3x_2}e^{w_1w_3(x_2+w_2y_2+\frac{w_2}{w_3}a)t}dx_2\\
&=\sum_{a=0}^{w_3d-1}(-1)^a\chi(a)\xi^{aw_1w_2}
(\sum_{k=0}^{\infty}E_{k,\chi,\xi^{w_2w_3}}(w_1y_1)\frac{(w_2y_3t)^k}{k!})\\
&\qquad\qquad\qquad\qquad\qquad\times
(\sum_{l=0}^{\infty}E_{l,\chi,\xi^{w_1w_3}}(w_2y_2+\frac{w_2}{w_3}a)\frac{(w_1y_3t)^l}{l!})\\
&=\sum_{n=0}^{\infty}(w_3^n\sum_{k=0}^{n}\binom{n}{k}E_{k,\chi,\xi^{w_2w_3}}(w_1y_1)
\sum_{a=0}^{w_3d-1}(-1)^a\chi(a)\xi^{aw_1w_2}\\
&\qquad\qquad\qquad\qquad\qquad\times E_{n-k,\chi,\xi^{w_1w_3}}
(w_2y_2+\frac{w_2}{w_3}a)w_1^{n-k}
w_2^k)\frac{t^n}{n!}.\\
&
\end{split}
\end{equation}

\noindent (a-2) Here we write $I(\Lambda_{23}^{2})$ in three
different ways:

\begin{equation}\label{N24}
\begin{split}
\noindent(1)~I(&\Lambda_{23}^{2})\\
&=\int_{X}\chi(x_1)\xi^{w_2w_3x_1}e^{w_2w_3(x_1+w_1y_1)t}dx_1\\
&\qquad\qquad\qquad\qquad\qquad\times\frac{\int_{X}\chi(x_2)\xi^{w_1w_3x_2}e^{w_1w_3x_2t}dx_2}{\int_{X}\xi^{dw_1w_2w_3x_4}e^{dw_1w_2w_3x_4t}dx_4}\\
&\qquad\qquad\qquad\qquad\qquad\qquad\qquad\times\frac{\int_{X}\chi(x_3)\xi^{w_1w_2x_3}e^{w_1w_2x_3t}dx_3}{\int_{X}\xi^{dw_1w_2w_3x_4}e^{dw_1w_2w_3x_4t}dx_4}\\
&=(\sum_{k=0}^{\infty}E_{k,\chi,\xi^{w_2w_3}}(w_1y_1)\frac{(w_2w_3t)^k}{k!})\\
&\qquad\qquad\qquad\qquad\times(\sum_{l=0}^{\infty}T_{l}(w_2d-1;\chi,\xi^{w_1w_3})\frac{(w_1w_3t)^l}{l!})\\
&\qquad\qquad\qquad\qquad\qquad\qquad\times
(\sum_{m=0}^{\infty}T_{m}(w_3d-1;\chi,\xi^{w_1w_2})\frac{(w_1w_2t)^m}{m!})
\\
\end{split}
\end{equation}

\begin{equation}\label{N25}
\begin{split}
=&\sum_{n=0}^{\infty}(\sum_{k+l+m=n}\binom{n}{k,l,m}E_{k,\chi,\xi^{w_2w_3}}(w_1y_1)T_l(w_2d-1;\chi,\xi^{w_1w_3})\\
&\qquad\qquad\qquad\qquad\qquad\quad\times
T_m(w_3d-1;\chi,\xi^{w_1w_2})w_{1}^{l+m}w_{2}^{k+m}w_{3}^{k+l})\frac{t^n}{n!}.
\end{split}
\end{equation}

\noindent (2) Invoking (\ref{N9}), (\ref{N24}) can also be written
as

\begin{equation}\label{N26}
\begin{split}
I(&\Lambda_{23}^{2})=\sum_{a=0}^{w_2d-1}(-1)^a\chi(a)\xi^{aw_1w_3}
\int_{X}\chi(x_1)\xi^{w_2w_3x_1}e^{w_2w_3(x_1+w_1y_1+\frac{w_1}{w_2}a)t}dx_1\\
&\qquad\qquad\qquad\qquad\qquad\qquad\qquad\qquad
\times\frac{\int_{X}\chi(x_3)\xi^{w_1w_2x_3}e^{w_1w_2x_3t}dx_3}{\int_{X}\xi^{dw_1w_2w_3x_4}e^{dw_1w_2w_3x_4t}dx_4}
\end{split}
\end{equation}
\begin{equation*}
\begin{split}
&\qquad\qquad\quad=\sum_{a=0}^{w_2d-1}(-1)^a\chi(a)\xi^{aw_1w_3}
(\sum_{k=0}^{\infty}E_{k,\chi,\xi^{w_2w_3}}(w_1y_1+\frac{w_1}{w_2}a)\frac{(w_2w_3t)^k}{k!})\\
&\qquad\qquad\qquad\qquad\qquad\qquad\qquad\qquad\times(\sum_{l=0}^{\infty}T_{l}(w_3d-1;\chi,\xi^{w_1w_2})\frac{(w_1w_2t)^l}{l!})\\
\end{split}
\end{equation*}

\begin{equation}\label{N27}
\begin{split}
&=\sum_{n=0}^{\infty}(w_2^n\sum_{k=0}^{n}\binom{n}{k}
\sum_{a=0}^{w_2d-1}(-1)^a\chi(a)\xi^{aw_1w_3}
E_{k,\chi,\xi^{w_2w_3}}(w_1y_1+\frac{w_1}{w_2}a)\\
&\qquad\qquad\qquad\qquad\qquad\qquad\times
T_{n-k}(w_3d-1;\chi,\xi^{w_1w_2})w_{1}^{n-k}w_3^{k})\frac{t^n}{n!}.
\end{split}
\end{equation}

\noindent (3) Invoking (\ref{N9}) once again, (\ref{N26}) can be
written as

\begin{equation*}
\begin{split}
I(&\Lambda_{23}^{2})\\
&=\sum_{a=0}^{w_2d-1}(-1)^a\chi(a)\xi^{aw_1w_3}
\sum_{b=0}^{w_3d-1}(-1)^b\chi(b)\xi^{bw_1w_2}\\
&\qquad\qquad\qquad\qquad\times\int_{X}\chi(x_1)\xi^{w_2w_3x_1}e^{w_2w_3(x_1+w_1y_1+\frac{w_1}{w_2}a+\frac{w_1}{w_3}b)t}dx_1\\
&=\sum_{a=0}^{w_2d-1}(-1)^a\chi(a)\xi^{aw_1w_3}
\sum_{b=0}^{w_3d-1}(-1)^b\chi(b)\xi^{bw_1w_2}\\
&\qquad\qquad\qquad\qquad\times\sum_{n=0}^{\infty}E_{n,\chi,\xi^{w_2w_3}}(w_1y_1+\frac{w_1}{w_2}a+\frac{w_1}{w_3}b)
\frac{(w_2w_3t)^n}{n!}
\end{split}
\end{equation*}

\begin{equation}\label{N28}
\begin{split}
=\sum_{n=0}^{\infty}((w_2w_3)^n\sum_{a=0}^{w_2d-1}&\sum_{b=0}^{w_3d-1}(-1)^{a+b}
\chi(ab)\xi^{w_1(aw_3+bw_2)}\\
&\qquad\qquad\times
E_{n,\chi,\xi^{w_2w_3}}(w_1y_1+\frac{w_1}{w_2}a+\frac{w_1}{w_3}b))
\frac{t^n}{n!}.
\end{split}
\end{equation}

\noindent (a-3)
\begin{equation*}
\begin{split}
I(\Lambda_{23}^{3})
&=\frac{\int_{X}\chi(x_1)\xi^{w_2w_3x_1}e^{w_2w_3x_1t}dx_1}{\int_{X}\xi^{dw_1w_2w_3x_4}e^{dw_1w_2w_3x_4t}dx_4}\\
&\qquad\qquad\qquad\times\frac{\int_{X}\chi(x_2)\xi^{w_1w_3x_2}e^{w_1w_3x_2t}dx_2}{\int_{X}\xi^{dw_1w_2w_3x_4}e^{dw_1w_2w_3x_4t}dx_4}\\
&\qquad\qquad\qquad\qquad\qquad\qquad\quad\times\frac{\int_{X}\chi(x_3)\xi^{w_1w_2x_3}e^{w_1w_2x_3t}dx_3}{\int_{X}\xi^{dw_1w_2w_3x_4}e^{dw_1w_2w_3x_4t}dx_4}\\
&=(\sum_{k=0}^{\infty}T_{k}(w_1d-1;\chi,\xi^{w_2w_3})\frac{(w_2w_3t)^k}{k!})\\
&\qquad\qquad\qquad\times(\sum_{l=0}^{\infty}T_{l}(w_2d-1;\chi,\xi^{w_1w_3})\frac{(w_1w_3t)^l}{l!})\\
&\qquad\qquad\qquad\qquad\qquad\times
(\sum_{m=0}^{\infty}T_{m}(w_3d-1;\chi,\xi^{w_1w_2})\frac{(w_1w_2t)^m}{m!})
\\
\end{split}
\end{equation*}

\begin{equation}\label{N29}
\begin{split}
=&\sum_{n=0}^{\infty}(\sum_{k+l+m=n}\binom{n}{k,l,m}T_k(w_1d-1;\chi,\xi^{w_2w_3})
T_l(w_2d-1;\chi,\xi^{w_1w_3})\qquad\\
&\qquad\qquad\qquad\qquad\qquad\quad\times
T_m(w_3d-1;\chi,\xi^{w_1w_2})
w_{1}^{l+m}w_{2}^{k+m}w_{3}^{k+l})\frac{t^n}{n!}.\\
\end{split}
\end{equation}

\noindent (b) For Type $\Lambda_{13}^{i}~($i=0,1,2,3$)$, we may
consider the analogous things to the ones in (a-0), (a-1), (a-2),
and (a-3). However, each of those can be obtained from the
corresponding ones in (a-0), (a-1), (a-2), and (a-3). Indeed, if we
substitute $w_2w_3$, $w_1w_3$, $w_1w_2$ respectively for $w_1$,
$w_2$, $w_3$ in (\ref{N11}), this amounts to replacing $t$ by
$w_1w_2w_3t$ in (\ref{N13}) and $\xi$ by $\xi^{w_1w_2w_3}$ in
(\ref{N13}). So, upon replacing $w_1$, $w_2$, $w_3$ respectively by
$w_2w_3$, $w_1w_3$, $w_1w_2$, dividing by $(w_1 w_2 w_3)^n$, and
replacing $\xi^{w_1w_2w_3}$ by $\xi$, in each of the expressions of
(\ref{N19}), (\ref{N22}), (\ref{N23}), (\ref{N25}),
(\ref{N27})-(\ref{N29}), we will get the corresponding symmetric
identities for Type $\Lambda_{13}^{i}$ ($i=0,1,2,3$).

\noindent (c-0)
\begin{equation*}
\begin{split}
I(&\Lambda_{12}^{0})\\
&=\int_{X}\chi(x_1)\xi^{w_1x_1}e^{w_1(x_1+w_2y)t}dx_1
\int_{X}\chi(x_2)\xi^{w_2x_2}e^{w_2(x_2+w_3y)t}dx_2\\
&\qquad\qquad\qquad\qquad\qquad\qquad\qquad\qquad
\times\int_{X}\chi(x_3)\xi^{w_3x_3}e^{w_3(x_3+w_1y)t}dx_3\\
&=(\sum_{k=0}^{\infty}\frac{E_{k,\chi,\xi^{w_1}}(w_2y)}{k!}(w_1t)^k)
(\sum_{l=0}^{\infty}\frac{E_{l,\chi,\xi^{w_2}}(w_3y)}{l!}(w_2t)^l)\\
&\qquad\qquad\qquad\qquad\qquad\qquad\qquad\qquad \times
(\sum_{m=0}^{\infty}\frac{E_{m,\chi,\xi^{w_3}}(w_1y)}{m!}(w_3t)^m)
\end{split}
\end{equation*}

\begin{equation}\label{N30}
\begin{split}
&=\sum_{n=0}^{\infty}(\sum_{k+l+m=n}\binom{n}{k,l,m}
E_{k,\chi,\xi^{w_1}}(w_2y) E_{l,\chi,\xi^{w_2}}(w_3y)\\
&\qquad\qquad\qquad\qquad\qquad\qquad\qquad\qquad\qquad \times
E_{m,\chi,\xi^{w_3}}(w_1y)w_{1}^{k}w_{2}^{l}w_{3}^{m}
)\frac{t^n}{n!}.
\end{split}
\end{equation}

\noindent (c-1)
\begin{equation*}
\begin{split}
I(&\Lambda_{12}^{1})\\
&=\frac{\int_{X}\chi(x_1)\xi^{w_1x_1}e^{w_1x_1t}dx_1}{\int_{X}\xi^{dw_1w_2z_3}e^{dw_1w_2z_3t}dz_3}
\times\frac{\int_{X}\chi(x_2)\xi^{w_2x_2}e^{w_2x_2t}dx_2}{\int_{X}\xi^{dw_2w_3z_1}e^{dw_2w_3z_1t}dz_1}\\
&\qquad\qquad\qquad\qquad\qquad\qquad\qquad\qquad
\times\frac{\int_{X}\chi(x_3)\xi^{w_3x_3}e^{w_3x_3t}dx_3}{\int_{X}\xi^{dw_3w_1z_2}e^{dw_3w_1z_2t}dz_2}\\
&=(\sum_{k=0}^{\infty}T_{k}(w_2d-1;\chi,\xi^{w_1})\frac{(w_1t)^k}{k!})
(\sum_{l=0}^{\infty}T_{l}(w_3d-1;\chi,\xi^{w_2})\frac{(w_2t)^l}{l!})\\
&\qquad\qquad\qquad\qquad\qquad\qquad\qquad\qquad\times
(\sum_{m=0}^{\infty}T_{m}(w_1d-1;\chi,\xi^{w_3})\frac{(w_3t)^m}{m!})
\\
\end{split}
\end{equation*}

\begin{equation}\label{N31}
\begin{split}
=&\sum_{n=0}^{\infty}(\sum_{k+l+m=n}\binom{n}{k,l,m}T_k(w_2d-1;\chi,\xi^{w_1})
T_l(w_3d-1;\chi,\xi^{w_2})\\
&\qquad\qquad\qquad\qquad\qquad\qquad\quad\times
T_m(w_1d-1;\chi,\xi^{w_3})w_{1}^{k}w_{2}^{l}w_{3}^{m})\frac{t^n}{n!}.\qquad\quad\\
&\\
\end{split}
\end{equation}

\section{Main theorems}

As we noted earlier in the last paragraph of Section 2, the various
types of quotients of $p$-adic fermionic integrals are invariant
under any permutation of $w_1$, $w_2$, $w_3$. So the corresponding
expressions in Section 3 are also invariant under any permutation of
$w_1$, $w_2$, $w_3$. Thus our results about identities of symmetry
will be immediate consequences of this observation.

However, not all permutations of an expression in Section 3 yield
distinct ones. In fact, as these expressions are obtained by
permuting $w_1$, $w_2$, $w_3$ in a single one labeled by them, they
can be viewed as a group in a natural manner and hence it is
isomorphic to a quotient of $S_3$. In particular, the number of
possible distinct expressions are 1,2,3, or 6. (a-0), (a-1(1)),
(a-1(2)), and (a-2(2)) give the full six identities of symmetry,
(a-2(1)) and (a-2(3)) yield three identities of symmetry, and (c-0)
and (c-1) give two identities of symmetry, while the expression in
(a-3) yields no identities of symmetry.

Here we will just consider the cases of Theorems \ref{T4} and
\ref{T8}, leaving the others as easy exercises for the reader. As
for the case of Theorem \ref{T4}, in addition to
(\ref{N42})-(\ref{N44}), we get the following three ones:

\begin{equation}\label{N32}
\begin{split}
\sum_{k+l+m=n}\binom{n}{k,l,m}E_{k,\chi,\xi^{w_2w_3}}(w_1y_1)&T_l(w_3d-1;\chi,\xi^{w_1w_2})\\
\times &T_m(w_2d-1;\chi,\xi^{w_1w_3})
w_{1}^{l+m}w_{3}^{k+m}w_{2}^{k+l},
\end{split}
\end{equation}
\begin{equation}\label{N33}
\begin{split}
\sum_{k+l+m=n}\binom{n}{k,l,m}E_{k,\chi,\xi^{w_1w_3}}(w_2y_1)&T_l(w_1d-1;\chi,\xi^{w_2w_3})\\
\times &T_m(w_3d-1;\chi,\xi^{w_1w_2})
w_{2}^{l+m}w_{1}^{k+m}w_{3}^{k+l},
\end{split}
\end{equation}
\begin{equation}\label{N34}
\begin{split}
\sum_{k+l+m=n}\binom{n}{k,l,m}E_{k,\chi,\xi^{w_1w_2}}(w_3y_1)&T_l(w_2d-1;\chi,\xi^{w_1w_3})\\
\times &T_m(w_1d-1;\chi,\xi^{w_2w_3})
w_{3}^{l+m}w_{2}^{k+m}w_{1}^{k+l}.
\end{split}
\end{equation}

\noindent But, by interchanging $l$ and $m$, we see that
(\ref{N32}), (\ref{N33}), and (\ref{N34}) are respectively equal to
(\ref{N42}), (\ref{N43}), and (\ref{N44}).

\noindent As to Theorem \ref{T8}, in addition to (\ref{N48}) and
(\ref{N49}), we have:

\begin{equation}\label{N35}
\begin{split}
\sum_{k+l+m=n}\binom{n}{k,l,m}T_k(w_2d-1;\chi,\xi^{w_1})&T_l(w_3d-1;\chi,\xi^{w_2})\\
\times & T_m(w_1d-1;\chi,\xi^{w_3}) w_{1}^{k}w_{2}^{l}w_{3}^{m},
\end{split}
\end{equation}

\begin{equation}\label{N36}
\begin{split}
\sum_{k+l+m=n}\binom{n}{k,l,m}T_k(w_3d-1;\chi,\xi^{w_2})&T_l(w_1d-1;\chi,\xi^{w_3})\\
\times & T_m(w_2d-1;\chi,\xi^{w_1}) w_{2}^{k}w_{3}^{l}w_{1}^{m},
\end{split}
\end{equation}

\begin{equation}\label{N37}
\begin{split}
\sum_{k+l+m=n}\binom{n}{k,l,m}T_k(w_3d-1;\chi,\xi^{w_1})&T_l(w_2d-1;\chi,\xi^{w_3})\\
\times & T_m(w_1d-1;\chi,\xi^{w_2}) w_{1}^{k}w_{3}^{l}w_{2}^{m},
\end{split}
\end{equation}

\begin{equation}\label{N38}
\begin{split}
\sum_{k+l+m=n}\binom{n}{k,l,m}T_k(w_2d-1;\chi,\xi^{w_3})&T_l(w_1d-1;\chi,\xi^{w_2})\\
\times & T_m(w_3d-1;\chi,\xi^{w_1}) w_{3}^{k}w_{2}^{l}w_{1}^{m}.
\end{split}
\end{equation}

\noindent However, (\ref{N35}) and (\ref{N36}) are equal to
(\ref{N48}), as we can see by applying the permutations
$k\rightarrow l$, $l\rightarrow m$, $m\rightarrow k$ for (\ref{N35})
and $k\rightarrow m$, $l \rightarrow k$, $m\rightarrow l$ for
(\ref{N36}). Similarly, we see that (\ref{N37}) and (\ref{N38}) are
equal to (\ref{N49}), by applying permutations $k\rightarrow l$,
$l\rightarrow m$, $m\rightarrow k$ for (\ref{N37}) and $k\rightarrow
m$, $l\rightarrow k$, $m\rightarrow l$ for (\ref{N38}).

\begin{theorem}\label{T1}
Let $w_1$, $w_2$, $w_3$ be any positive integers. Then we have:

\begin{equation*}
\begin{split}
\sum_{k+l+m=n}\binom{n}{k,l,m}E_{k,\chi,\xi^{w_2w_3}}(w_1y_1)&E_{l,\chi,\xi^{w_1w_3}}(w_2y_2)\\
\times&E_{m,\chi,\xi^{w_1w_2}}(w_3y_3)
w_{1}^{l+m}w_{2}^{k+m}w_{3}^{k+l}\\
=\sum_{k+l+m=n}\binom{n}{k,l,m}E_{k,\chi,\xi^{w_2w_3}}(w_1y_1)&E_{l,\chi,\xi^{w_1w_2}}(w_3y_2)\\
\times&E_{m,\chi,\xi^{w_1w_3}}(w_2y_3)
w_{1}^{l+m}w_{3}^{k+m}w_{2}^{k+l}\\
\end{split}
\end{equation*}

\begin{equation}\label{N39}
\begin{split}
=\sum_{k+l+m=n}\binom{n}{k,l,m}E_{k,\chi,\xi^{w_1w_3}}(w_2y_1)&E_{l,\chi,\xi^{w_2w_3}}(w_1y_2)\\
\times&E_{m,\chi,\xi^{w_1w_2}}(w_3y_3)
w_{2}^{l+m}w_{1}^{k+m}w_{3}^{k+l}\\
=\sum_{k+l+m=n}\binom{n}{k,l,m}E_{k,\chi,\xi^{w_1w_3}}(w_2y_1)&E_{l,\chi,\xi^{w_1w_2}}(w_3y_2)\\
\times&E_{m,\chi,\xi^{w_2w_3}}(w_1y_3)
w_{2}^{l+m}w_{3}^{k+m}w_{1}^{k+l}\\
=\sum_{k+l+m=n}\binom{n}{k,l,m}E_{k,\chi,\xi^{w_1w_2}}(w_3y_1)&E_{l,\chi,\xi^{w_2w_3}}(w_1y_2)\\
\times&E_{m,\chi,\xi^{w_1w_3}}(w_2y_3)
w_{3}^{l+m}w_{1}^{k+m}w_{2}^{k+l}\\
=\sum_{k+l+m=n}\binom{n}{k,l,m}E_{k,\chi,\xi^{w_1w_2}}(w_3y_1)&E_{l,\chi,\xi^{w_1w_3}}(w_2y_2)\\
\times&E_{m,\chi,\xi^{w_2w_3}}(w_1y_3)
w_{3}^{l+m}w_{2}^{k+m}w_{1}^{k+l}.
\end{split}
\end{equation}
\end{theorem}

\begin{theorem}\label{T2}
Let $w_1$, $w_2$, $w_3$ be any odd positive integers. Then we have:
\begin{equation}\label{N40}
\begin{split}
\sum_{k+l+m=n}\binom{n}{k,l,m}E_{k,\chi,\xi^{w_2w_3}}(w_1y_1)&E_{l,\chi,\xi^{w_1w_3}}(w_2y_2)\\
\times&T_m(w_3d-1;\chi,\xi^{w_1w_2})
w_{1}^{l+m}w_{2}^{k+m}w_{3}^{k+l}\\
=\sum_{k+l+m=n}\binom{n}{k,l,m}E_{k,\chi,\xi^{w_2w_3}}(w_1y_1)&E_{l,\chi,\xi^{w_1w_2}}(w_3y_2)\\
\times&T_m(w_2d-1;\chi,\xi^{w_1w_3})
w_{1}^{l+m}w_{3}^{k+m}w_{2}^{k+l}\\
=\sum_{k+l+m=n}\binom{n}{k,l,m}E_{k,\chi,\xi^{w_1w_3}}(w_2y_1)&E_{l,\chi,\xi^{w_2w_3}}(w_1y_2)\\
\times&T_m(w_3d-1;\chi,\xi^{w_1w_2})
w_{2}^{l+m}w_{1}^{k+m}w_{3}^{k+l}\\
=\sum_{k+l+m=n}\binom{n}{k,l,m}E_{k,\chi,\xi^{w_1w_3}}(w_2y_1)&E_{l,\chi,\xi^{w_1w_2}}(w_3y_2)\\
\times&T_m(w_1d-1;\chi,\xi^{w_2w_3})
w_{2}^{l+m}w_{3}^{k+m}w_{1}^{k+l}\\
=\sum_{k+l+m=n}\binom{n}{k,l,m}E_{k,\chi,\xi^{w_1w_2}}(w_3y_1)&E_{l,\chi,\xi^{w_1w_3}}(w_2y_2)\\
\times&T_m(w_1d-1;\chi,\xi^{w_2w_3})
w_{3}^{l+m}w_{2}^{k+m}w_{1}^{k+l}\\
=\sum_{k+l+m=n}\binom{n}{k,l,m}E_{k,\chi,\xi^{w_1w_2}}(w_3y_1)&E_{l,\chi,\xi^{w_2w_3}}(w_1y_2)\\
\times&T_m(w_2d-1;\chi,\xi^{w_1w_3})
w_{3}^{l+m}w_{1}^{k+m}w_{2}^{k+l}.
\end{split}
\end{equation}
\end{theorem}

\begin{theorem}\label{T3}
Let $w_1$, $w_2$, $w_3$ be any odd positive integers. Then we have:
\begin{equation}\label{N41}
\begin{split}
w_1^n\sum_{k=0}^{n}\binom{n}{k}E_{k,\chi,\xi^{w_1w_2}}(w_3y_1)&\sum_{a=0}^{w_1d-1}(-1)^a\chi(a)\xi^{aw_2w_3}\\
\times& E_{n-k,\chi,\xi^{w_1w_3}}(w_2y_2+\frac{w_2}{w_1}a)w_3^{n-k}w_2^k\\
=w_1^n\sum_{k=0}^{n}\binom{n}{k}E_{k,\chi,\xi^{w_1w_3}}(w_2y_1)&\sum_{a=0}^{w_1d-1}(-1)^a\chi(a)\xi^{aw_2w_3}\\
\times& E_{n-k,\chi,\xi^{w_1w_2}}(w_3y_2+\frac{w_3}{w_1}a)w_2^{n-k}w_3^k\\
=w_2^n\sum_{k=0}^{n}\binom{n}{k}E_{k,\chi,\xi^{w_1w_2}}(w_3y_1)&\sum_{a=0}^{w_2d-1}(-1)^a\chi(a)\xi^{aw_1w_3}\\
\times& E_{n-k,\chi,\xi^{w_2w_3}}(w_1y_2+\frac{w_1}{w_2}a)w_3^{n-k}w_1^k\\
=w_2^n\sum_{k=0}^{n}\binom{n}{k}E_{k,\chi,\xi^{w_2w_3}}(w_1y_1)&\sum_{a=0}^{w_2d-1}(-1)^a\chi(a)\xi^{aw_1w_3}\\
\times& E_{n-k,\chi,\xi^{w_1w_2}}(w_3y_2+\frac{w_3}{w_2}a)w_1^{n-k}w_3^k\\
=w_3^n\sum_{k=0}^{n}\binom{n}{k}E_{k,\chi,\xi^{w_1w_3}}(w_2y_1)&\sum_{a=0}^{w_3d-1}(-1)^a\chi(a)\xi^{aw_1w_2}\\
\times & E_{n-k,\chi,\xi^{w_2w_3}}(w_1y_2+\frac{w_1}{w_3}a)w_2^{n-k}w_1^k\\
=w_3^n\sum_{k=0}^{n}\binom{n}{k}E_{k,\chi,\xi^{w_2w_3}}(w_1y_1)&\sum_{a=0}^{w_3d-1}(-1)^a\chi(a)\xi^{aw_1w_2}\\
\times &
E_{n-k,\chi,\xi^{w_1w_3}}(w_2y_2+\frac{w_2}{w_3}a)w_1^{n-k}w_2^k.
\end{split}
\end{equation}
\end{theorem}

\begin{theorem}\label{T4}
Let $w_1$, $w_2$, $w_3$ be any odd positive integers. Then we have
the following three symmetries in $w_1$, $w_2$, $w_3$:
\begin{equation}\label{N42}
\begin{split}
\sum_{k+l+m=n}\binom{n}{k,l,m}E_{k,\chi,\xi^{w_2w_3}}&(w_1y_1)T_l(w_2d-1;\chi,\xi^{w_1w_3})\\
\times & T_m(w_3d-1;\chi,\xi^{w_1w_2})
w_{1}^{l+m}w_{2}^{k+m}w_{3}^{k+l}
\end{split}
\end{equation}
\begin{equation}\label{N43}
\begin{split}
=\sum_{k+l+m=n}\binom{n}{k,l,m}E_{k,\chi,\xi^{w_1w_3}}&(w_2y_1)T_l(w_3d-1;\chi,\xi^{w_1w_2})\\
\times & T_m(w_1d-1;\chi,\xi^{w_2w_3})
w_{2}^{l+m}w_{3}^{k+m}w_{1}^{k+l}\\
\end{split}
\end{equation}
\begin{equation}\label{N44}
\begin{split}
=\sum_{k+l+m=n}\binom{n}{k,l,m}E_{k,\chi,\xi^{w_1w_2}}&(w_3y_1)T_l(w_1d-1;\chi,\xi^{w_2w_3})\\
\times & T_m(w_2d-1;\chi,\xi^{w_1w_3})
w_{3}^{l+m}w_{1}^{k+m}w_{2}^{k+l}.
\end{split}
\end{equation}
\end{theorem}

\begin{theorem}\label{T5}
Let $w_1$, $w_2$, $w_3$ be any odd positive integers. Then we have:
\begin{equation}\label{N45}
\begin{split}
w_1^n\sum_{k=0}^{n}\binom{n}{k}\sum_{a=0}^{w_1d-1}(-1)^a\chi(a)\xi^{w_2w_3}
&E_{k,\chi,\xi^{w_1w_3}}(w_2y_1+\frac{w_2}{w_1}a)\\
\times &T_{n-k}(w_3d-1;\chi,\xi^{w_1w_2})w_2^{n-k}w_3^k\\
=w_1^n\sum_{k=0}^{n}\binom{n}{k}\sum_{a=0}^{w_1d-1}(-1)^a\chi(a)\xi^{w_2w_3}
&E_{k,\chi,\xi^{w_1w_2}}(w_3y_1+\frac{w_3}{w_1}a)\\
\times &T_{n-k}(w_2d-1;\chi,\xi^{w_1w_3})w_3^{n-k}w_2^k\\
=w_2^n\sum_{k=0}^{n}\binom{n}{k}\sum_{a=0}^{w_2d-1}(-1)^a\chi(a)\xi^{w_1w_3}
&E_{k,\chi,\xi^{w_2w_3}}(w_1y_1+\frac{w_1}{w_2}a)\\
\times &T_{n-k}(w_3d-1;\chi,\xi^{w_1w_2})w_1^{n-k}w_3^k\\
=w_2^n\sum_{k=0}^{n}\binom{n}{k}\sum_{a=0}^{w_2d-1}(-1)^a\chi(a)\xi^{w_1w_3}
&E_{k,\chi,\xi^{w_1w_2}}(w_3y_1+\frac{w_3}{w_2}a)\\
\times &T_{n-k}(w_1d-1;\chi,\xi^{w_2w_3})w_3^{n-k}w_1^k\\
=w_3^n\sum_{k=0}^{n}\binom{n}{k}\sum_{a=0}^{w_3d-1}(-1)^a\chi(a)\xi^{w_1w_2}
&E_{k,\chi,\xi^{w_2w_3}}(w_1y_1+\frac{w_1}{w_3}a)\\
\times &T_{n-k}(w_2d-1;\chi,\xi^{w_1w_3})w_1^{n-k}w_2^k\\
=w_3^n\sum_{k=0}^{n}\binom{n}{k}\sum_{a=0}^{w_3d-1}(-1)^a\chi(a)\xi^{w_1w_2}
&E_{k,\chi,\xi^{w_1w_3}}(w_2y_1+\frac{w_2}{w_3}a)\\
\times &T_{n-k}(w_1d-1;\chi,\xi^{w_2w_3})w_2^{n-k}w_1^k.
\end{split}
\end{equation}
\end{theorem}

\begin{theorem}\label{T6}
Let $w_1$, $w_2$, $w_3$ be any odd positive integers. Then we have
the following three symmetries in $w_1$, $w_2$, $w_3$:
\begin{equation}\label{N46}
\begin{split}
&(w_1w_2)^n\sum_{a=0}^{w_1d-1}\sum_{b=0}^{w_2d-1}(-1)^{a+b}\chi(ab)\xi^{w_3(aw_2+bw_1)}
E_{n,\chi,\xi^{w_1w_2}}(w_3y_1+\frac{w_3}{w_1}a+\frac{w_3}{w_2}b)\\
=&(w_2w_3)^n\sum_{a=0}^{w_2d-1}\sum_{b=0}^{w_3d-1}(-1)^{a+b}\chi(ab)\xi^{w_1(aw_3+bw_2)}
E_{n,\chi,\xi^{w_2w_3}}(w_1y_1+\frac{w_1}{w_2}a+\frac{w_1}{w_3}b)\\
=&(w_3w_1)^n\sum_{a=0}^{w_3d-1}\sum_{b=0}^{w_1d-1}(-1)^{a+b}\chi(ab)\xi^{w_2(aw_1+bw_3)}
E_{n,\chi,\xi^{w_1w_3}}(w_2y_1+\frac{w_2}{w_3}a+\frac{w_2}{w_1}b).
\end{split}
\end{equation}
\end{theorem}

\begin{theorem}\label{T7}
Let $w_1$, $w_2$, $w_3$ be any positive integers. Then we have the
following two symmetries in $w_1$, $w_2$, $w_3$:
\begin{equation}\label{N47}
\begin{split}
&\sum_{k+l+m=n}\binom{n}{k,l,m}E_{k,\chi,\xi^{w_3}}(w_1y)E_{l,\chi,\xi^{w_1}}(w_2y)
E_{m,\chi,\xi^{w_2}}(w_3y)
w_{3}^{k}w_{1}^{l}w_{2}^{m}\\
=&\sum_{k+l+m=n}\binom{n}{k,l,m}E_{k,\chi,\xi^{w_2}}(w_1y)E_{l,\chi,\xi^{w_1}}(w_3y)
E_{m,\chi,\xi^{w_3}}(w_2y) w_{2}^{k}w_{1}^{l}w_{3}^{m}.
\end{split}
\end{equation}
\end{theorem}

\begin{theorem}\label{T8}
Let $w_1$, $w_2$, $w_3$ be any odd positive integers. Then we have
the following two symmetries in $w_1$, $w_2$, $w_3$:
\begin{equation}\label{N48}
\begin{split}
\sum_{k+l+m=n}\binom{n}{k,l,m}T_k(w_1d-1;\chi,\xi^{w_3})&T_l(w_2d-1;\chi,\xi^{w_1})\\
\times &T_m(w_3d-1;\chi,\xi^{w_2}) w_{3}^{k}w_{1}^{l}w_{2}^{m}
\end{split}
\end{equation}
\begin{equation}\label{N49}
\begin{split}
\sum_{k+l+m=n}\binom{n}{k,l,m}T_k(w_1d-1;\chi,\xi^{w_2})&T_l(w_3d-1;\chi,\xi^{w_1})\\
\times &T_m(w_2d-1;\chi,\xi^{w_3}) w_{2}^{k}w_{1}^{l}w_{3}^{m}.
\end{split}
\end{equation}
\end{theorem}


\end{document}